\documentclass[preprint,12pt]{elsarticle}%

\usepackage{etex}
\usepackage[english]{babel}
\usepackage[utf8x]{inputenc}
\usepackage{amsmath}
\usepackage{amsfonts}
\usepackage{stmaryrd}
\usepackage{amssymb}
\usepackage{subfigure}
\usepackage{graphicx}
\usepackage{multicol}
\usepackage[colorinlistoftodos]{todonotes}
\usepackage{fullpage}

\makeatletter
\def\ps@pprintTitle{%
 \let\@oddhead\@empty
 \let\@evenhead\@empty
 \def\@oddfoot{}%
 \let\@evenfoot\@oddfoot}
\makeatother

\newcommand{\calT}{\mathcal{T}}

\newcommand{\bfE}{\mathbf{E}}
\newcommand{\bfH}{\mathbf{H}}

\newcommand{\bfn}{\mathbf{n}}
\newcommand{\bft}{\mathbf{t}}
\newcommand{\bfv}{\mathbf{v}}
\newcommand{\bfV}{\mathbf{V}}

\newcommand{\bfJ}{\mathbf{J}}

\newcommand{\bfr}{\mathbf{r}}
\newcommand{\iu}{\mathrm{i}}
\newcommand{\calF}{\mathcal{F}}
\newcommand{\jump}[1]{\llbracket #1 \rrbracket}
\newcommand{\average}[1]{\left\{ #1 \right\}}
\newcommand{\ds}{\displaystyle}

\begin{document}
\begin{frontmatter}

\title{A hybridizable discontinuous Galerkin method for solving nonlocal optical response models}

\author[lab1]{Liang Li\corref{cor1}} 
\ead{plum\_liliang@uestc.edu.cn, plum.liliang@gmail.com}
\author[lab2]{St\'ephane Lanteri} 
\ead{stephane.lanteri@inria.fr}
\author[lab3,lab4]{N. Asger Mortensen} 
\ead{asger@mail.aps.org}
\author[lab3,lab4]{Martijn Wubs} 
\ead{mwubs@fotonik.dtu.dk}

\address[lab1]{School of Mathematical Sciences, University of
               Electronic Science and Technology of China, 611731,
               Chengdu, P.R. China}

\address[lab2]{INRIA, 2004 Route des Lucioles, BP 93  
               06902 Sophia Antipolis Cedex, France}

\address[lab3]{Department of Photonics Engineering, Technical University of Denmark,
{\O}rsteds Plads 343, DK-2800 Kgs. Lyngby, Denmark.}

\address[lab4]{Center for Nanostructured Graphene, Technical University of Denmark,
{\O}rsteds Plads 343, DK-2800 Kgs. Lyngby, Denmark} 

\begin{abstract}
We propose Hybridizable Discontinuous Galerkin (HDG) methods for solving the frequency-domain Maxwell's equations coupled to the Nonlocal Hydrodynamic Drude (NHD) and Generalized Nonlocal Optical Response (GNOR) models, which are employed to describe the optical properties of nano-plasmonic scatterers and waveguides. Brief derivations for both the NHD model and the GNOR model are presented. The formulations of the HDG method are given, in which we introduce two hybrid variables living only on the skeleton of the mesh. The local field solutions are expressed in terms of the hybrid variables in each element. Two conservativity conditions are globally enforced to make the problem solvable and to guarantee the continuity of the tangential component of the electric field and the normal component of the current density. Numerical results show that the proposed HDG methods converge at optimal rate. We benchmark our implementation and demonstrate that the HDG method has the potential to solve complex nanophotonic problems.

\end{abstract}


\end{frontmatter}

\section{Introduction}
\label{sec:intro}
Nanophotonics is the active research field field concerned with the study of
interactions between nanometer scale structures/media and light, including near-infrared, visible, and ultraviolet light. It bridges the micro and the macro worlds, and there are many connections between theoretical studies and feasible engineering. The many fascinating (potential) applications include invisibility cloaking, nano antennas, metamaterials, novel biological detection and treatment technologies, as well as new storage media~\cite{Maier}.  

All of the above applications of nanophotonics require elaborate control of  the propagation of light waves. In order to do so, appropriate mathematical models are needed to predict the behavior of light-matter interactions. Metals are interesting for nanophotonics because they can both enhance and confine optical fields, making plasmonics of interest to emerging quantum technologies~\cite{Tame,Bozhevolnyi,Fitzgerald}. This is enabled by the existence of Surface Plasmons (SPs). SPs  are coherent oscillations that exist as evanescent waves at both sides of the interface between any two materials where the real part of the dielectric function changes sign across the interface. The typical example is a metal-dielectric interface, such as a metal sheet in air~\cite{SaridChallener}. Maxwell's equations can be employed to model the macroscale electromagnetic waves and armed with classical electrodynamics there are numerous approaches ranging from classical electrodynamics to ab initio treatments~\cite{Gallinet2015,Varas2016}. Ab initio techniques can be used to simulate the microscopic dynamics on the atomic scale, but with ab initio methods one can only deal with systems with  up to about ten thousand atoms~\cite{Varas2016}, thus calling for semiclassical treatments~\cite{MortensenGNOR,ToscanoEtal2015} or more effective inclusions of quantum phenomena into classical electrodynamics~\cite{Luo2013,Yan2015,Christensen2016,Zhu2016}.

If one models the interaction of light with metallic nanostructures classically or semiclassically, then this calls for appropriate modelling of the material response as described for example by the Drude model~\cite{Drude1900, DresselScheffler2006}, the Nonlocal Hydrodynamic Drude (NHD) model~\cite{Bloch1933a,RazaEtal2015JPCM}, or the Generalized
Nonlocal Optical Response (GNOR) theory~\cite{MortensenGNOR}, all in combination with and coupled to Maxwell's equations. Except for some highly symmetric geometries, analytical solutions to the resulting systems of  differential equations are not available. Thus, numerical treatment of  these systems of PDEs is an important aspect of nanophotonics research.   Numerical experiments help to find promising systems and geometries before  real fabrication, to obtain optimized parameters, to visualize field distributions, to investigate the dominant contribution to a phenomenon, to explain experimental observations, and so on~\cite{Busch2011}.

Several  numerical methods exist for computing the solution of Maxwell's equations~\cite{Gallinet2015}. For time-dependent problems, the Finite-Difference Time-Domain (FDTD) algorithm is the most popular method~\cite{Taflove2005} among physicists and engineers. More recently, the Discontinuous Galerkin Time-Domain (DGTD) method has drawn a lot of attention because of several appealing features, for example, easy adaptation to complex geometries and material composition, high-order accuracy, and natural parallelism~\cite{HesthavenWarburton2002}. For time-harmonic problems, the Finite Element Method (FEM) is most widely used for the solution of Maxwell's equations. In very recent years, the Hybridizable  Discontinuous Galerkin (HDG) method appears as a promising numerical method  for time-harmonic problems because it inherits nearly all the advantages of the DG methods while leading to a computational complexity similar to FEM~\cite{CGL2009, NPC2011, LLP2014, HLLH2016}.

Currently, FDTD  (for time-dependent problems) and FEM (for time-harmonic problems) methods are still the methods most commonly adopted for the simulation of light-matter interactions. Most often, commercial simulation
software (such as Lumerical FDTD\footnote{https://www.lumerical.com/}
and Comsol Multiphysics\footnote{http://www.comsol.com/}) is used for that  purpose. However, these  methods and  computer codes do not always offer the required capabilities for addressing accurately and efficiently the complexity of the physical phenomena underlying nanometer scale light-matter interactions. In the academic community, also the DGTD method has recently been considered in this context~\cite{Busch2011,Viquerat2015,HuangYQ2016}. In Ref.~\cite{Schmittetal2016}, some numerical results are presented for the NHD model using the DGTD method. In the present paper we are employing the HDG method to solve the frequency-domain NHD and GNOR models. The development of accurate and efficient numerical  methods for computational  nanophotonics is expected to be a long-lasting demand, both because new models are regularly proposed that require innovative numerical methods, and because there is demand for more accurate and faster simulation methods for existing models.   

This paper introduces a HDG method for the solution the NHD and GNOR models. The rest of the paper is organized as follows. In section~\ref{sec:models}, we briefly introduce mathematical aspects both of the NHD model and of the GNOR model. HDG formulations are given in section~\ref{sec:hdg}. Numerical results are presented in section~\ref{sec:tests} to show the effectiveness of high-order HDG methods for solving problems in nanophotonics. We draw conclusions in section~\ref{sec:con}.

%---------------------------------------------------------------------
\section{Physics problem: nonlocal optical response by nanoparticles}
\label{sec:models}
\begin{figure}
\centering
\includegraphics[scale=0.5]{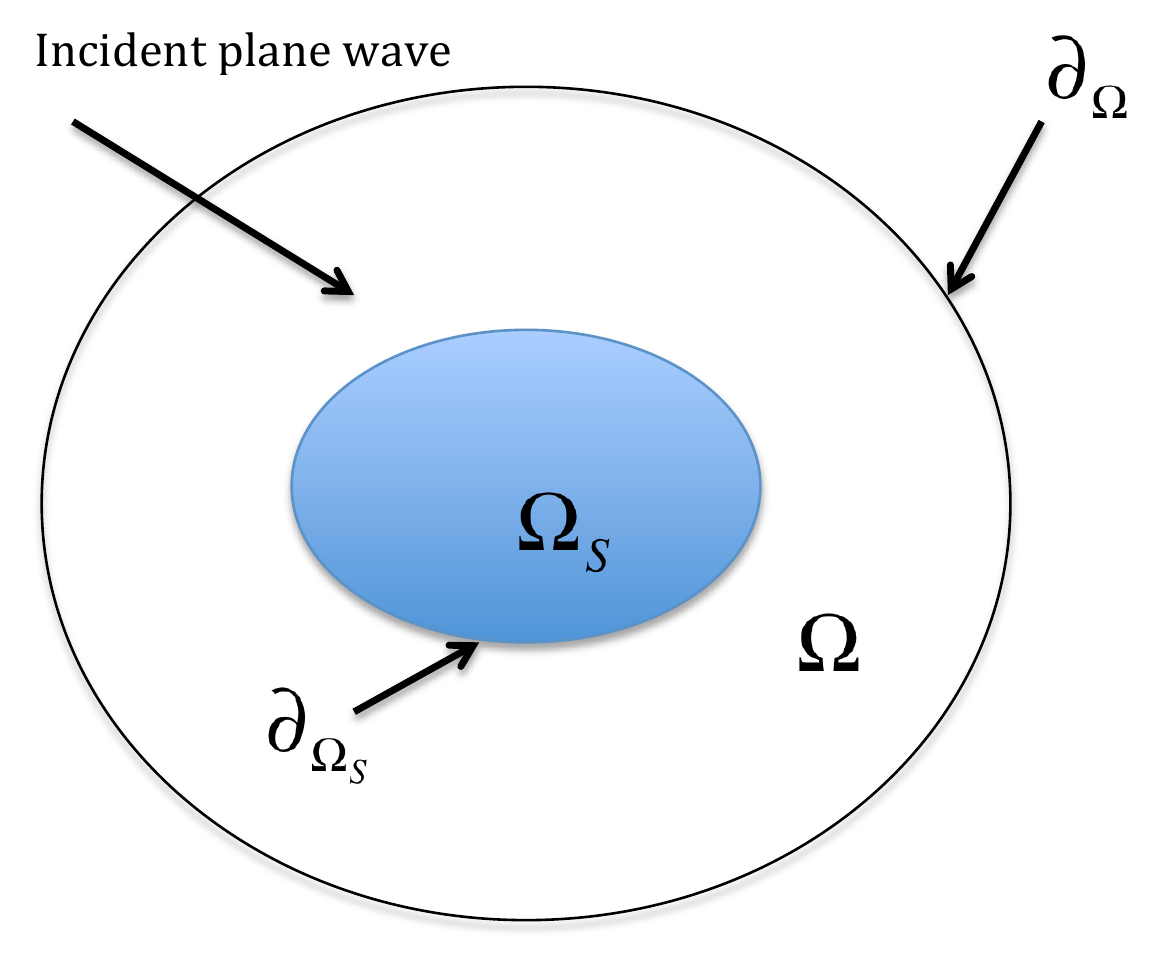}
\caption{Sketch of the incident electromagnetic wave illuminating the scatterer $\Omega_S$ that has a subwavelength size and is surrounded by free space. $\Omega_S$ is usually filled with metals, such as gold, silver or sodium. An artificial absorbing boundary $\partial\Omega$ is introduced to make a computational domain $\Omega$.}\label{fig:struct}
\end{figure}

The problem considered is shown in Figure~\ref{fig:struct} where the nanometer-size metal $\Omega_S$ is illuminated by an incident plane wave of light. The infinite scattering domain is truncated as a finite computational domain $\Omega$ by employing an artificial absorbing boundary condition, which is designed to absorb outgoing waves. 

\subsection{Nonlocal hydrodynamic Drude model}
There are a number of theories for the modeling of the light-matter interactions which are used under different settings. In this subsection, we briefly introduce the NHD model. The incoming light is described as a propagating electromagnetic wave that satisfies Maxwell's equations. Without external charge and current, Maxwell's equations of macroscopic electromagnetism for non-magnetic materials can be written as
\begin{equation}
\label{eq:max}
\left\{
\begin{array}{l}
\nabla\times \bfH =   \varepsilon_0\varepsilon_{\text{loc}}\partial_t\bfE+\bfJ,\\
\nabla\times\bfE = -\mu_0\partial_t\mathbf{H},
\end{array}
\right.
\end{equation}
where $\mathbf{H}$ and $\bfE$ are respectively the magnetic and electric fields, $ \varepsilon_0$ is the permittivity constant, $\mu_0$ is the permeability constant, $\varepsilon_{\text{loc}}=\varepsilon_{\infty}+\varepsilon_{\text{inter}}$ is introduced to account for the local response, and $\bfJ$ is the nonlocal hydrodynamic polarization current density which is due to the nonlocal material on the plasmonic scatterers~\cite{HiremathEtal2012}. In this paper, we will for simplicity set $\varepsilon_{\text{inter}}=0$ and $\varepsilon_{\infty}=1$, thereby focusing solely on the free-electron response to light. Equations \eqref{eq:max} need to be completed to solve electromagnetic fields $\bfE$ and $\mathbf{H}$ because of the unknown polarization current density $\bfJ$. The models that we will consider in this paper differ only in the assumed dynamics of the polarization current density, which we will now discuss in more detail. 

The polarization current density $\bfJ$ due to the motion of the free-electron gas can be written as
\begin{equation}
\label{eq:jHD}
\bfJ = -en\bfv,
\end{equation}
where $e$ is the charge of the electron, $n$ is the density of the electron gas (a scalar field), and $\bfv$ is its hydrodynamic velocity (a vector field). Within the hydrodynamic model, the dynamics of the velocity field is given by~\cite{Schmittetal2016,RazaEtal2015}
\begin{equation}
\label{eq:hd}
m_e(\partial_t+\bfv\cdot\nabla)\bfv=-e(\bfE+\bfv\times\mathbf{B}) - m_e\gamma\bfv -\nabla\left(\frac{\delta g[n]}{\delta n}\right),
\end{equation}
where $m_e$ is the mass of an electron, $-e(\bfE+\bfv\times\mathbf{B})$ is the Lorentz force with $\mathbf{B}$ being the magnetic flux density, $\gamma$ is a damping constant, $g[n]$ is an energy functional of the fluid, and the term $\nabla\left(\frac{\delta g[n]}{\delta n}\right)$ denotes the quantum pressure. Complementary to Eq.~\ref{eq:hd}, the dynamics of the free-electron density is given by 
\begin{equation}
\partial_t n+\nabla\cdot(n\bfv)=0,
\end{equation}
which is the well-known continuity relation that relates the velocity $\bfv$ and the density $n$.

The hydrodynamic dynamics described by Eq.~(\ref{eq:hd}) is obviously nonlinear in ${\bf v}$, but in the following we only consider the linear response of the electron gas on external fields. One can write a perturbation expansion ${\bf v} \simeq {\bf v}_{0} + {\bf v}_{1}$ and similarly for the electric and magnetic fields and for the density. Since in the absence of an external field ${\bf v} = {\bf v}_{0} = {\bf 0}$, both the nonlinear term $\bfv\cdot\nabla\bfv$ and the  magnetic induction field $\mathbf{B}$ disappear due to the linearization~\cite{ToscanoEtal2015}. If we furthermore assume the energy functional to be of the Thomas-Fermi form, then we obtain for the linearized quantum pressure 
\begin{equation}\label{eq:qpressure}
-\nabla\left(\frac{\delta g[n]}{\delta n}\right)= -m_e\beta^2\frac{1}{n_0}\nabla n,
\end{equation}
where $\beta^2=\frac{3}{5}v_F^2$ with $v_F$ being the Fermi velocity. The zero-order (i.e. equilibrium) density $n_0$ is constant within the plasmonic medium~\cite{ToscanoEtal2015}. Here in Eq.~\eqref{eq:qpressure} and below, we write $n$ for the linearized density ${n}_{1}$ and similarly we will from now on simply write ${\bf v}$ for the linearized velocity ${\bf v}_{1}$. As a result, we obtain the linearized hydrodynamic equation~\cite{RazaEtal2015JPCM,Schmittetal2016}
\begin{equation}
\label{eq:hds} 
m_e\partial_t\bfv=-e\bfE-m_e\gamma\bfv-m_e\beta^2\frac{1}{n_0}\nabla n,
\end{equation}
as well as the linearized continuity relation
\begin{equation}
\label{eq:nv}
\partial_t n = -n_0 \nabla\cdot \bfv.
\end{equation}
Inserting Eqs.~\eqref{eq:jHD} (linearized as $\bfJ=-en_0\bfv$) and \eqref{eq:nv} into \eqref{eq:hds}, and taking the time-derivative $\partial_t$, we obtain
\begin{equation}
\label{eq:hdT} 
\partial_{tt} \bfJ+\gamma\partial_t\bfJ-\beta^2\nabla(\nabla\cdot\bfJ)-\omega_p^2\varepsilon_0\partial_t\bfE=0,
\end{equation}
where $\omega_p$ is the plasma frequency with $\omega_p^2=n_0e^2/(m_e\varepsilon_0)$. By Fourier transformation we replace $\partial_t$ with $-\iu\omega$, where $\iu$ is the imaginary unit and $\omega$ is the angular frequency, and obtain the frequency-domain relation between polarization current density and the electric field within the hydrodynamic model as
\begin{equation}\label{eq:hdFre}
\omega(\omega+\iu\gamma)\bfJ+\beta^2\nabla(\nabla\cdot\bfJ)=\iu\omega\omega^2_p\varepsilon_0\bfE.
\end{equation}
This equation describes electron-field interaction within the plasmonic nanostructure $\Omega_S$. We will neglect spill-out of electrons outside the classical geometric surface of the structure, which for our purposes is a good assumption for noble metals such as silver and gold~\cite{ToscanoEtal2015}. Mathematically, this is arranged by imposing a hard-wall condition  on the boundary $\partial\Omega_S$, namely $\bfn\cdot\bfJ=0$ on $\partial\Omega_S$~\cite{Jewsbury:1981a,Yan:2013a}.

%------------------------------------------------------------

\subsection{General nonlocal optical response model}
We also briefly present the mathematical derivation of the central equations of the GNOR model, based on Ref.~\cite{MortensenGNOR}. In the GNOR model, also diffusion of the electron gas is taken into consideration. Let the density $n(\bfr,t)=n_0+n_1(\bfr,t)$, where the last term is the induced density variation caused by a non-vanishing electric field $\bfE$, which we assume sufficiently small that $n_1\ll n_0$ holds. Instead of~\eqref{eq:nv}, we now consider the linearized convection-diffusion equation~\cite{MortensenGNOR}
\begin{equation}
\label{eq:cd}
\partial_t (-e_1) n_1 = D\nabla^2(-e)n_1-\nabla\cdot[(-e)n_0\bfv]=-\nabla\cdot\bfJ,
\end{equation}
where $D$ is the diffusion constant for the charge-carrier diffusion. Then the current density is given by Fick's law
\begin{equation}
\label{eq:current}
\bfJ = (-e)n_0\bfv-D\nabla(-e)n_1.
\end{equation} 
Multiplying \eqref{eq:hds} by the charge of the electron $-e$, the equilibrium density $n_0$ and taking the time-derivative we have
\begin{equation}
\label{eq:hddv}
m_e(\partial_t+\gamma)\partial_t[(-e)n_0\bfv]=n_0e^2\partial_t\bfE-m_e\beta^2\nabla[\partial_t(-e)n_1].
\end{equation}
Dividing \eqref{eq:hddv} by $m_e$ and combining with Fick's law~\eqref{eq:current} results in
\begin{equation}
\label{eq:hddj}
(\partial_t+\gamma)\{\partial_t\bfJ+D\nabla[\partial_t(-e)n_1]\}=\frac{n_0e^2}{m_e}\partial_t\bfE-\beta^2\nabla[\partial_t(-e)n_1].
\end{equation}
From the convection-diffusion equation~\eqref{eq:cd}, we have
\begin{equation}
\label{eq:hdd}
(\partial_t+\gamma)[\partial_t\bfJ+D\nabla(\nabla\cdot\bfJ)]=\frac{n_0e^2}{m_e}\partial_t\bfE-\beta^2\nabla(\nabla\cdot\bfJ).
\end{equation}
Like what we did for~\eqref{eq:hdT}, transforming \eqref{eq:hdd} to the frequency domain gives 
\begin{equation}\label{eq:gnor}
\omega(\omega+\iu\gamma)\bfJ+[\beta^2+D(\gamma-\iu\omega)]\nabla(\nabla\cdot\bfJ)=\iu\omega\omega^2_p\varepsilon_0\bfE.
\end{equation}
The physical predictions obtained by the GNOR and NHD models often differ substantially, as illustrated below. However, from a computational point of view  the GNOR model only differs by the replacement $\beta^2\rightarrow\beta^2+D(\gamma-\iu\omega)$ in the frequency domain, whereby the nonlocal hydrodynamic parameter acquires an often non-negligible imaginary part. In the GNOR model we have the same additional boundary condition $\bfn\cdot\bfJ=0$ on $\partial\Omega_S$ as in the NHD model.
%------------------------------------------------------------

\subsection{Specification to 2D TM mode}
Now we can couple Maxwell's equation \eqref{eq:max} with \eqref{eq:hdFre} for the NHD model, or similarly with \eqref{eq:gnor} for the GNOR model. We will compute light extinction by infinitely long nanowires. We take the wire axes along the $z$-direction and consider TM-polarized incident light, i.e. polarized in the $(x,y)$-plane. In this 2D setting, we can define $\bfE = (E_x, E_y)^T$ to be a vector and $H = H_z$ a scalar function. Coupling the time-harmonic Maxwell's equations and hydrodynamic Drude model \eqref{eq:hdFre}, we have in 2D
\begin{equation}\label{eq:hyd}
\left\{
\begin{aligned}
\nabla\times H & =  -\iu\omega\varepsilon_0\bfE +\bfJ, \text{ in } \Omega, \\
\nabla\times\bfE & = \iu\omega\mu_0 H, \text{ in } \Omega,\\
\nabla(\nabla\cdot\bfJ) + \frac{\omega(\omega+\iu\gamma)}{\beta^2}\bfJ & =\frac{\iu\omega\omega^2_p\varepsilon_0}{\beta^2}\bfE, \text{ in } \Omega_S.
\end{aligned}
\right.
\end{equation}
If the Silver-M{\"u}ller boundary condition (first-order absorbing boundary condition) \cite{Stupfel1994} is applied on the boundary $\partial\Omega$ of the computational domain, then we have the boundary conditions
\begin{equation}
\label{eq:bc}
\left\{
\begin{aligned}
\bfn\times\bfE-H & = \bfn\times\bfE^{\text{inc}}-H^{\text{inc}} = g^{\text{inc}}, \text{ on } \partial\Omega, \\
\bfn\cdot\bfJ & =0, \text{ on } \partial\Omega_S,
\end{aligned}
\right.
\end{equation}
where $\bfE^{\text{inc}}$ and $H^{\text{inc}}$ stand for the electromagnetic fields of the incoming light.
%------------------------------------------------------------
\section{HDG formulations of nonlocal optical response models }\label{sec:hdg}

\subsection{The promise of hybridizable DG methods}\label{Sec:promiseHDG}
In the Introduction some properties and advantages of DG and HDG methods were briefly mentioned, which we here explain in more detail. The classic DG method is seldomly employed for solving stationary problems, because it duplicates degrees of freedom (DOFs) on every internal edge. Thus the number of globally coupled DOFs is much greater than the number of DOFs required by conforming finite element methods for the same accuracy. Consequently, DG methods are expensive in terms of both CPU time and memory consumption. Hybridization of DG methods~\cite{CGL2009} is devoted to addressing this issue while at the same time keeping all the advantages of DG methods. HDG methods introduce additional hybrid variables on the edges of the elements. Then we define the numerical traces arising from partial integration in the DG formulations through the hybrid variables. We can thus define the local (element-wise) solutions by hybrid variables. Conservativity conditions are imposed on numerical traces to ensure the continuity of the tangential component of the electric field and the normal component of the current density and to make the problem solvable. As a result, HDG methods produce a linear system in terms of the DOFs of the additional hybrid variables only. In this way, the number of globally coupled DOFs is greatly reduced as compared to the classic DG method. In a recent study~\cite{Yakovlev:2016}, the authors showed that HDG methods outperform FEM in many cases.

\subsection{Computational concepts and notations}
\begin{figure}
\centering
\includegraphics[scale=0.5]{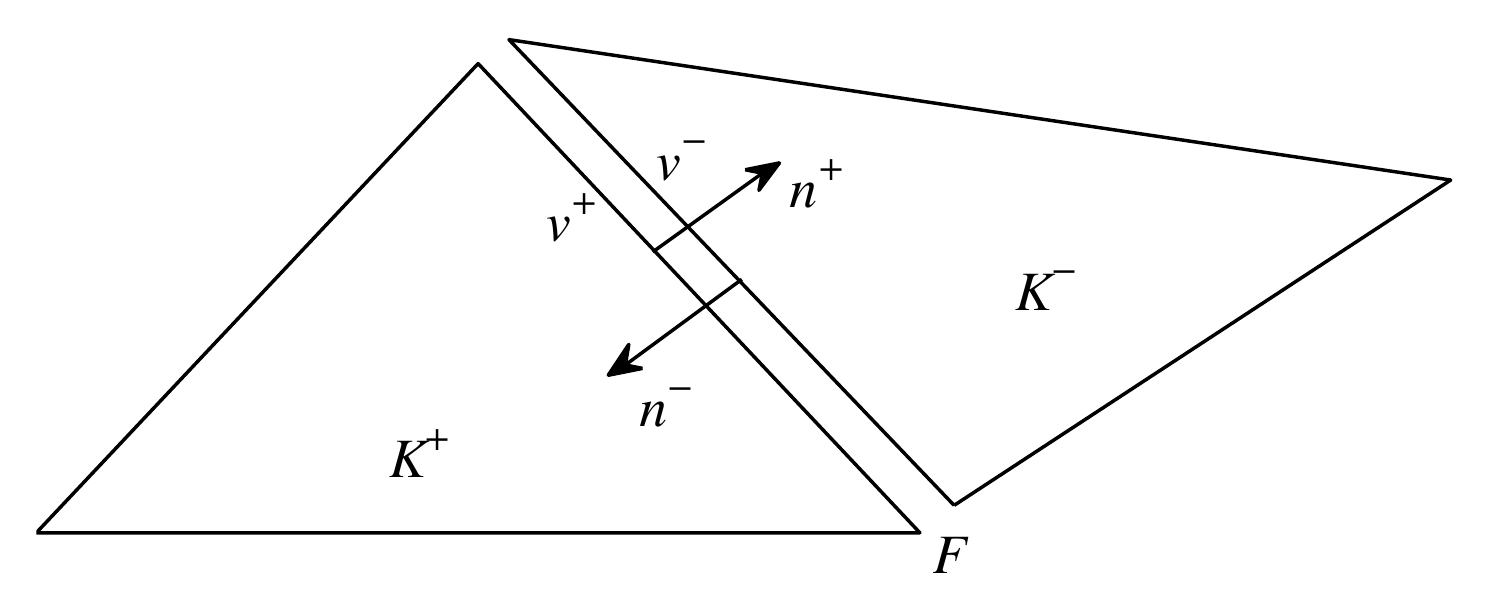}
\caption{Two neighboring discretization elements (here: triangles) within the computational domain. An edge $F$ is shared by two elements $K^+$ and $K^-$. The outward normal vectors $n^+$ and $n^-$ point in opposite directions. A characteristic property of the DG method is that computed functions are allowed to be discontinuous across $F$ (hence the ``D'' in DG). For example, for a function $v$, be it a scalar or a vector, its value on $F$ from $K^+$ is $v^+$, while its value on $F$ from $K^-$ is $v^-$, and these $v^+$ and $v^-$ are not necessarily equal. By contrast, the hybrid variables in the HDG method {\em are} single-valued on~$F$.}\label{fig:interEdge}
\end{figure}

In order to give a clear presentation of the HDG method, here we introduce some computational concepts and notations. We divide the computational domain $\Omega$ into triangle elements. The union of all the triangles is denoted by $\calT_h$. By $\calF_h$ we denote the union of all edges of $\calT_h$. Furthermore, $\calF^I_h$ stands for the union of all the edges associated with the nanostructure. For an edge associated with two elements $F =\overline{K^+}\cap\overline{K^-}\in\calF_h$, let $(\bfv^\pm, v^\pm)$ be the \emph{traces} of $(\bfv, v)$ on $F$ from the  interior of $K^{\pm}$, see Fig.~\ref{fig:interEdge}, where we use the term \emph{trace} to denote the restriction of a function on the boundaries of the elements~\cite{Arnold2002}. Note that from now on $\bfv$ is used to describe a general vector function instead of velocity.
On  every face,  we define  \emph{mean}  (\emph{average}) \emph{values} $\average \cdot$ and \emph{jumps} $\jump \cdot$ as 
\begin{displaymath}
  \left\{
    \begin{aligned}
      \average \bfv_F  & = \dfrac{1}{2} (\bfv^+ + \bfv^-),
      \\
      \average v_F &  = \dfrac{1}{2}(v^+ + v^-),
      \\
      \jump{\bfn\times\bfv}_F &  = \bfn^+\times\bfv^+ +
      \bfn^-\times\bfv^-, \\
       \jump{\bfn\cdot\bfv}_F &  = \bfn^+\cdot\bfv^+ +
      \bfn^-\cdot\bfv^-,
      \\
      \jump {v\bft}_F &  = v^+ \bft^+ + v^- \bft^-,
    \end{aligned}
  \right. 
\end{displaymath}
where $\bfn^{\pm}$ denotes the outward unit norm vector to $K^\pm$ and $\bft^{\pm}$ denotes the unit tangent vectors to the boundaries $\partial    K^{\pm}$ such that $\bft^+\times\bfn^+=1$    and $\bft^-\times\bfn^-=1$. For the boundary edges, either on $\partial\Omega$ or on $\partial\Omega_S$, these expressions are modified as
\begin{displaymath}
  \left\{
    \begin{aligned}
      \average \bfv_F &  = \bfv^+,
      \\
      \average v_F &  = v^+,
      \\
      \jump{\bfn\times\bfv}_F &  = \bfn^+ \times\bfv^+, \\
       \jump{\bfn\cdot\bfv}_F & = \bfn^+ \cdot\bfv^+,
      \\
      \jump{v\bft}_F & = v^+ \bft^+.
    \end{aligned}
  \right .
\end{displaymath}
Let $\mathbb{P}_p(D)$ denote the space of polynomial functions of degree at most $p$ on a  domain $D$. For any element $K\in\calT_h$, let $V^p(K)$ be the space $\mathbb{P}_p(K)$ and $\bfV^p(K)$ the space $(\mathbb{P}_p(K))^2$. The discontinuous finite element spaces are then defined by
\begin{equation}\label{eq:dgSpace}
  \begin{aligned}
    V^p_h & = \left\{ v \in L^{2}(\Omega) \ \lvert \ v|_K\in
    V^p(K), \ \forall K\in\calT_h\right\},
    \\
    \bfV^p_h & = \left\{ \bfv \in (L^2(\Omega))^2 \ \lvert \
    \bfv|_K\in \bfV^p(K), \ \forall K\in\calT_h \right\},
  \end{aligned}
\end{equation}
where $L^2(\Omega)$ is the space of square integrable functions on the domain $\Omega$. We also introduce a traced finite element space
$$
\label{SpaceM}
M^p_h = \left\{ \eta \in L^2(\calF_h) \ \lvert
  \ \eta|_F \in\mathbb{P}_p(F),
  \ \forall F\in\calF_h  \right\}.
$$
Note that $M^p_h$ consists of functions which are continuous on an edge, but discontinuous at its ends. The restrictions of $V^p_h$, $\bfV^p_h$ and $M^p_h$ in $\Omega_S$ are denoted by $\widetilde{V}^p_h$, $\widetilde{\bfV}^p_h$ and $\widetilde{M}^p_h$.
For two vectorial functions $\mathbf{u}$ and $\bfv$ in $(L^2(D))^2$, we introduce the inner product $(\mathbf{u},\bfv)_{D}   =   \ds\int_D\mathbf{u}\cdot\overline{\bfv}  \,  {\rm d}x$,  where $\overline{\cdot}$ denotes the complex conjugation. Likewise for scalar functions $u$ and $v$ in $L^2(D)$, the inner product is defined as $(u,v)_D = \ds\int_D u  \overline v \,{\rm d}x$ provided $D$ is a  domain in $\mathbb{R}^2$. Finally we define the edge overlap $\langle  u,v\rangle_F =\ds\int_F u \overline v \, {\rm d}s$, where $F$ is a specific edge. Accordingly, we can define the total edge overlap for the whole triangulation or for relevant subsets of edges. Important cases are
\begin{displaymath}
    \langle\cdot,\cdot\rangle_{\calF_h}  =
    \sum_{F\in\calF_h}\langle\cdot,\cdot\rangle_F, \qquad
     \langle\cdot,\cdot\rangle_{\partial\Omega}  =
    \sum_{F\in\calF_h\cap\partial\Omega}\langle\cdot,\cdot\rangle_F,
     \qquad
    \langle\cdot,\cdot\rangle_{\calF^I_h}  =
    \sum_{F\in\calF^I_h}\langle\cdot,\cdot\rangle_F.
\end{displaymath}
denoting, respectively, the total edge overlap on the computational domain, the cumulative edge overlap on the absorbing boundary of the computational domain, and finally the cumulative edge overlap on the nanostructure. 
%---------------------------------------------------------------------%

\subsection{DG formulation of the coupled electrodynamical  equations }
\label{subsec:dg}
We begin the construction of a DG implementation of the hydrodynamic Drude model by rewriting the coupled electrodynamical equations~\eqref{eq:hyd} into a system of first-order equations
\begin{equation}
\label{eq:MaxHyd}
\left\{
\begin{aligned}
\iu\omega \varepsilon_0\bfE + \nabla\times H -\bfJ &= 0 \qquad \text{ in } \Omega, \\
\iu\omega\mu_0 H - \nabla\times\bfE & =0 \qquad \text{ in } \Omega, \\
\nabla q +\frac{ \gamma-\iu\omega}{\beta^2}\bfJ - \frac{\omega_p^2\varepsilon_0}{\beta^2}\bfE & = 0 \qquad \text{ in } \Omega_S, \\
\iu\omega q - \nabla\cdot\bfJ & = 0 \qquad \text{ in } \Omega_S,
\end{aligned}
\right.
\end{equation}
where we introduced the scalar function 
$q = (\iu\omega)^{-1}\nabla\cdot\bfJ $ which coincides with a scaled charge density.
In general, a DG method seeks an approximate solution $(\bfE_h, H_h, \bfJ_h, q_h)$ in the space $\bfV^{p}_{h}\times V^{p}_{h}\times\widetilde{\bfV}^{p}_{h}\times \widetilde{V}^{p}_{h}$ that for each element $K$ (in our case: for each discretization triangle) satisfies~\cite{ElBouajaji2013}
\begin{equation}
\label{eq:MaxHydInt}
\left\{
\begin{aligned}
(\iu\omega \varepsilon_0\bfE_h, \bfv)_{K} + (\nabla\times H_h, \bfv)_{K} - (\bfJ_h, \bfv)_{K} & = 0 \qquad \forall \bfv\in \bfV^{p}(K), \\
(\iu\omega\mu_0 H_h, v)_{K} - (\nabla\times\bfE_h, v)_{K} & =0 \qquad \forall v\in V^{p}(K), \\
(\nabla q_h, \bfv)_{K} + \Big(\frac{\gamma-\iu\omega}{\beta^2}\bfJ_h, \bfv\Big)_{K} - \Big(\frac{\omega_p^2\varepsilon_0}{\beta^2}\bfE_h, \bfv\Big)_{K} & = 0 \qquad \forall \bfv\in \widetilde{\bfV}^{p}(K), \\
(\iu\omega q_h, v)_{K} - (\nabla\cdot\bfJ_h, v)_{K} & = 0 \qquad \forall v\in \widetilde{V}^{p}(K).
\end{aligned}
\right.
\end{equation}
The application of appropriate Green's formulas to this system of equations leads to terms on the element boundaries~\cite{Arnold2002}. These boundary terms are the keys to connect the elements, since the elements themselves are independent due to the nature of the discontinuous finite elements spaces of Eq.~\eqref{eq:dgSpace}. In a DG method, one replaces the boundary terms by so-called \emph{numerical traces} $\hat{\bfE}_{h}, \hat{H}_h, \hat{\bfJ}_{h}$ and $\hat{q}_{h}$~\cite{CGL2009,LLP2013}, which are also known as `numerical fluxes' in the literature~\cite{Busch2011}. 
These numerical traces are defined as
\begin{equation}\label{eq:dgNumericalTrace}
\left\{
\begin{aligned}
\hat{H}_h & =\average {H_h} + \alpha_{E}\jump{\bfn\times\bfE_h}, \\
\bfn\times\hat{\bfE}_h & = \average{\bfn\times\bfE_h}+\alpha_{H}\jump{H_h},\\
\hat{q}_h & = \average{q_h}+\alpha_J\jump{\bfn\cdot\bfJ_h},\\
\bfn\cdot\hat{\bfJ}_h & =\average{\bfn\cdot\bfJ}+\alpha_q\jump{q_h}.
\end{aligned}
\right.
\end{equation}
In these definitions there is still freedom to choose values for the $\alpha$ parameters, and this corresponds to different DG schemes: by setting $\alpha_{E}=\alpha_{H}=\alpha_{J}=\alpha_{q}=0$, one obtains the \emph{centered flux} DG scheme. With $\alpha_{E}=\alpha_{H}=\alpha_{J}=\alpha_{q}=1$, one obtains the \emph{upwind flux} DG scheme~\cite{ElBouajaji2013}. For more validated DG schemes, we refer the interested readers to Ref.~\cite{Arnold2002}.
Having defined the numerical traces, we finally form a global system of linear equations involving all the DOFs on all the elements
\begin{equation}
\label{eq:MaxHydIntGreen}
\left\{
\begin{aligned}
(\iu\omega \varepsilon_0\bfE_h, \bfv)_{K} + (H_h, \nabla\times \bfv)_{K} -  \langle\hat{H}_{h},\bfn\times\bfv\rangle_{\partial K} - (\bfJ_h, \bfv)_{K} & = 0 \quad  \forall \bfv\in \bfV^{p}(K), \\
(\iu\omega\mu_0 H_h, v)_{K} - (\bfE_h, \nabla\times v)_{K} - \langle\bfn\times\hat{\bfE}_{h}, v\rangle_{\partial K} & =0 \quad \forall v\in V^{p}(K), \\
-(q_h, \nabla\cdot\bfv)_{K} + \langle \hat{q}_{h}, \bfn\cdot\bfv\rangle_{\partial K} +  \Big(\frac{\gamma-\iu\omega}{\beta^2}\bfJ_h, \bfv\Big)_{K} - \Big(\frac{\omega_p^2\varepsilon_0}{\beta^2}\bfE_h, \bfv\Big)_{K} & = 0 \quad \forall \bfv\in \widetilde\bfV^{p}(K), \\
(\iu\omega q_h, v)_{K} + (\bfJ_h, \nabla v)_{K} -\langle\bfn\cdot\hat{\bfJ}_{h}, v\rangle_{\partial K} & = 0 \quad \forall v\in\widetilde V^{p}(K),
\end{aligned}
\right.
\end{equation}
which are coupled equations that are valid whatever DG scheme is adopted.

\subsection{Hybridizable DG implementation of the electrodynamical  equations }
\label{subsec:hdg}
In Sec.~\ref{Sec:promiseHDG} we mentioned that hybridized DG methods have advantages as compared to the classic DG schemes, and here we discuss the hybridized approach in more detail. Unlike in the above classic DG formulations where the numerical traces directly couple the values from the elements on both sides of the edges, in a HDG formulation the numerical traces are defined through hybrid variables. Introducing two hybrid variables $\lambda_h$ and $\eta_h$ which live only on the boundaries of the elements, we define the numerical traces by
\begin{equation}
\label{eq:numTrac}
\left\{
\begin{aligned}
\hat{H}_h & = \lambda_h, \\
\hat{\bfE}_h & = \bfE_h + \tau_{\lambda}(\lambda_h-H_h)\bft, \\
\hat{q}_h & = \eta_h, \\
\hat{\bfJ}_h & = \bfJ + \tau_{\eta}(q_h-\eta_h)\bfn,
\end{aligned}
\right.
\end{equation}
where $\tau_{\lambda}$ and $\tau_{\eta}$ are two stabilization parameters.
Replacing the numerical traces in \eqref{eq:MaxHydIntGreen} with the expressions in \eqref{eq:numTrac} and applying Green's formulas to the first and fourth equations in \eqref{eq:MaxHydIntGreen}, we obtain the local formulation of the HDG method as
\begin{equation}
\label{eq:MaxHydIntGreen2}
\left\{
\begin{aligned}
(\iu\omega \varepsilon_0\bfE_h, \bfv)_{K} + (H_h, \nabla\times \bfv)_{K} - \langle\lambda_{h},\bfn\times\bfv\rangle_{\partial K} - (\bfJ_h, \bfv)_{K} & = 0, \ \forall \bfv\in \bfV^{p}(K), \\
(\iu\omega\mu_0 H_h, v)_{K} - (\nabla\times\bfE_h, v)_{K} + \langle\tau_{\lambda}(H_h-\lambda_h), v\rangle_{\partial K} & =0, \ \forall v\in V^{p}(K), \\
-(q_h, \nabla\cdot\bfv)_{K} + \langle\eta_{h}, \bfn\cdot\bfv\rangle_{\partial K} +  \Big(\frac{\gamma-\iu\omega}{\beta^2}\bfJ_h, \bfv\Big)_{K} - \Big(\frac{\omega_p^2\varepsilon_0}{\beta^2}\bfE_h, \bfv\Big)_{K} & = 0, \ \forall \bfv\in\widetilde \bfV^{p}(K), \\
(\iu\omega q_h, v)_{K} - (\nabla\cdot\bfJ_h, v)_{K} -\langle\tau_{\eta}(q_h-\eta_h), v\rangle_{\partial K} & = 0, \ \forall v\in\widetilde V^{p}(K)  .
\end{aligned}
\right.
\end{equation}
One can solve the local fields element by element once the solutions for $\lambda_h$ and $\eta_h$ are obtained.
In order to make the problem solvable, we need to employ global conditions
\begin{equation}\label{eq:globalRaw}
\left\{
\begin{aligned}
\langle\jump{\bfn\times\hat{\bfE}_h}, v\rangle_{\calF_h} - \langle\lambda_h, v\rangle_{\partial\Omega} & = \langle g^{\text{inc},v}\rangle_{\partial\Omega}, \ \forall v\in M^p_h, \\
\langle\jump{\bfn\cdot\hat{\bfJ}_h}, v\rangle_{\calF_h^I} & = 0, \ \forall v\in\widetilde M^p_h.
\end{aligned}
\right.
\end{equation}
The first relation in~\eqref{eq:globalRaw} weakly enforces the continuity condition for the tangential component of the electric field across any edges, and also takes into account the Silver-M{\" u}ller absorbing boundary condition. The other global condition in Eq.~\eqref{eq:globalRaw} weakly enforces the continuity condition for the normal component of the current density across any edges. The additional boundary condition on the surface of the nanostructure is implicitly contained in this relation.

Substituting $\hat{\bfE}_h$ and $\hat{\bfJ}_h$ in \eqref{eq:globalRaw} with the definitions in \eqref{eq:numTrac}, we arrive at the global reduced system of equations
\begin{equation}
\label{eq:global}
\left\{
\begin{aligned}
\langle\bfn\times\bfE_h-\tau_{\lambda}(\lambda_h-H_h), v\rangle_{\calF_h} - \langle\lambda_h, v\rangle_{\partial\Omega} & =  \langle g^{\text{inc}},v\rangle_{\partial\Omega}, \ \forall v\in M^p_h, \\
\langle\bfn\cdot\bfJ_h+\tau_\eta(q_h-\eta_h), v\rangle_{\calF_h^I} & = 0, \ \forall v\in\widetilde M^p_h.
\end{aligned}
\right.
\end{equation}
Note that we used the fact that $\bfn\times\bft=-1$ in \eqref{eq:global}. The two relations in Eq.~\eqref{eq:global} are not independent. They are coupled through the local solutions of $\bfE$, $H$, $\bfJ$ and $q$ of the local equations~\eqref{eq:MaxHydIntGreen2}.

\noindent\textbf{Remark I.} The proposed HDG formulation for the global system \eqref{eq:global} is naturally consistent with the boundary conditions, both on the artificial boundary and on the medium boundary. 

\noindent\textbf{Remark II.} Globally, we only need to solve Eq.~\eqref{eq:global}, in which the fields $\bfE_h$, $H_h$, $\bfJ_h$ and $q_h$ are replaced by the solutions in terms of $\lambda_h$ and $\eta_h$ from the local problems~\eqref{eq:MaxHydIntGreen2}. So the global DOFs are  associated with $\lambda_h$ in the whole computational domain, while they are associate with $\eta_h$ only within in the material medium. The discretization leads to a system of linear equations
\begin{equation}
\label{eq:linSys}
A\begin{bmatrix}
\underline{\lambda}_h \\[5pt]
\underline{\eta}_h
\end{bmatrix}=\begin{bmatrix}
\underline{g}^\text{inc}_h \\[5pt]
0
\end{bmatrix},
\end{equation}
where $\underline{\lambda}_h$ and $\underline{\eta}_h$ are vectors accounting for the degrees of freedom of the hybrid variables $\lambda$ and $\eta$ respectively, and the coefficient matrix $A$ is large and sparse.

%---------------------------------------------------------------------
\section{Numerical results}
\label{sec:tests}
In this section we present numerical results to validate the proposed HDG formulations. All HDG methods have been implemented in Fortran 90. All our tests are performed on a Macbook with a 1.3 GHz Inter Core i5 CPU and 4 GB memory. We employ the multifrontal sparse direct solver MUMPS~\cite{AmestoyEtal2000} to solve the discretized systems of linear equations. 

In HDG methods, we calculate the total fields $\bfE^\text{tot}$ and $H^\text{tot}$. The scattered fields are then calculated by subtracting the incident field from the total fields. We use HDG-$\mathbb{P}_p$ to denote the HDG method with interpolation order $p$. Here we choose fixed values $\tau_{\lambda}= \tau_{\eta}=1$ for the stabilization parameters. Different choices are discussed in Ref.~\cite{Gopalakrishnan2015}.  

\subsection{Convergence study: Wave propagation in a cavity}
\begin{figure}
\centering
\includegraphics[scale=0.45]{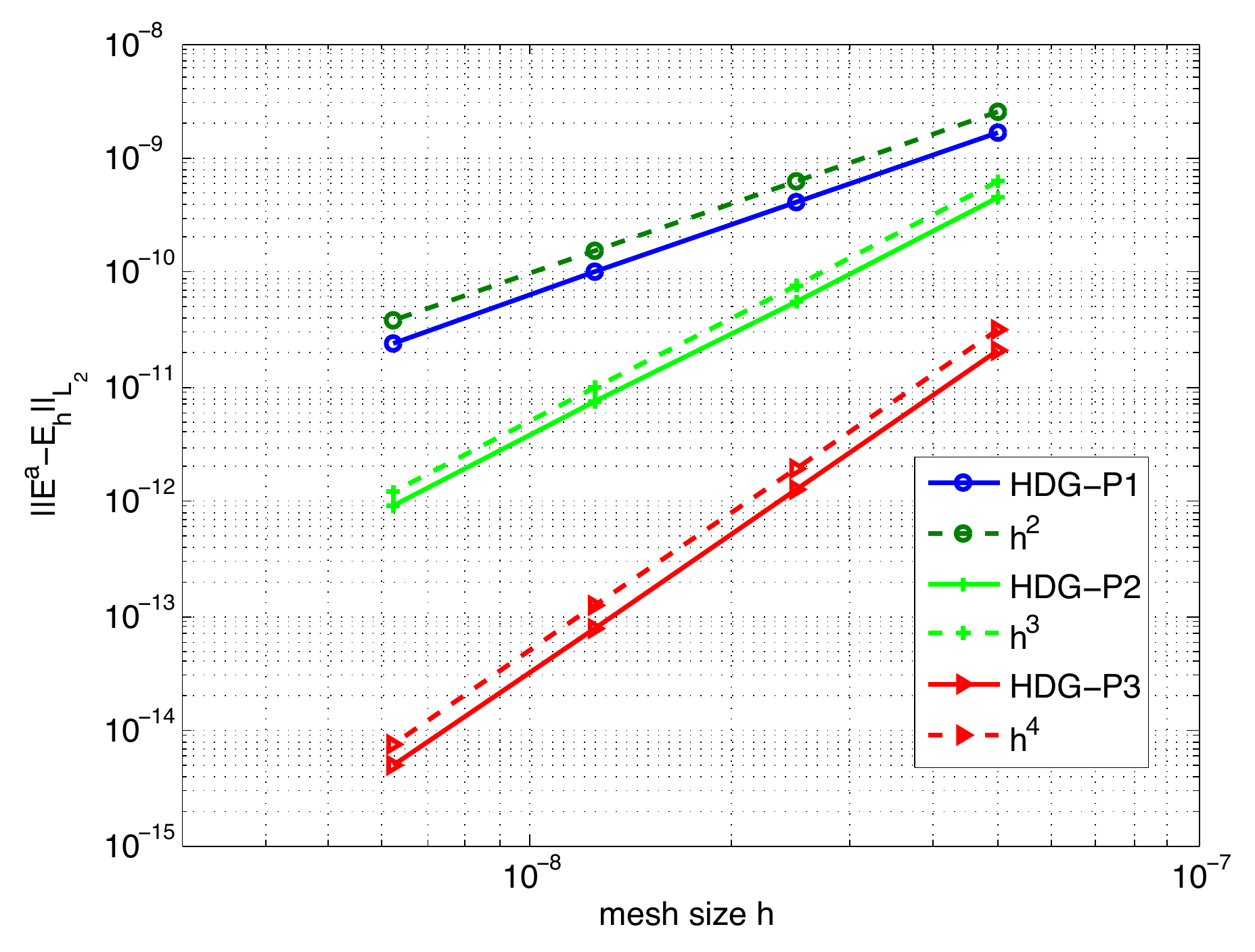}
\caption{Convergence history of the proposed HDG method for the cavity problem.}\label{fig:convRect}
\end{figure}

While elsewhere in this article we focus on nanowire structures, here we first study the convergence of our method by considering wave propagation in a cavity. This cavity is assumed to be a square domain $\Omega_{\Box}=\{(x,y)\in[0,L]\times[0,L]\}$ with the PEC boundary condition and hard-wall condition
$$
\bfn \times \bfE = 0, \text{ and } \bfn\cdot\bfJ=0, \text { on } \partial\Omega_{\Box}.
$$
This test case can be viewed as the frequency-domain version of the first test case in~\cite{Schmittetal2016}. The simplicity is achieved by introducing artificial current density and electric field, such that the analytical solutions coincide with Maxwell's equations and with the hydrodynamic equation
\begin{equation}
\label{eq:MaxHydArt}
\left\{
\begin{aligned}
\iu k\bfE + \nabla\times H & = \bfJ -\bfJ^{a}, \\
\iu k H - \nabla\times\bfE & =0, \\
\nabla q -\frac{ \iu\omega}{\beta^2}\bfJ &  = -\frac{\gamma}{\beta^2}\bfJ+\frac{\gamma}{\beta^2}\bfJ^{a} + \frac{\omega_p^2\varepsilon_0}{\beta^2}\bfE -\frac{\omega_p^2\varepsilon_0}{\beta^2}\bfE^{a}, \\
\iu\omega q - \nabla\cdot\bfJ &  = 0,
\end{aligned}
\right.
\end{equation}
where $k=\frac{\omega}{c}$ is the wave number, with $c$ being the light speed. We make this modification to unify the scale of the electric and magnetic fields. The artificial terms $\bfJ^a$ and $\bfE^a$ are also the analytical solution to this equation~\eqref{eq:MaxHydArt}:
\begin{equation}\label{eq:anaSol}
\begin{aligned}
\bfE^a & = \frac{\sqrt{2}}{2}\iu\begin{bmatrix}
-\cos(\frac{\sqrt{2}}{2}kx)\sin(\frac{\sqrt{2}}{2}ky)\\[5pt]
\sin(\frac{\sqrt{2}}{2}kx)\cos(\frac{\sqrt{2}}{2}ky)
\end{bmatrix},\\%%%
\bfJ^a & = -\frac{\sqrt{2}\mu_0 k\beta^2}{2\omega}\begin{bmatrix}
\sin(\frac{\sqrt{2}}{2}kx)\cos(\frac{\sqrt{2}}{2}ky)\\[5pt]
\cos(\frac{\sqrt{2}}{2}kx)\sin(\frac{\sqrt{2}}{2}ky)
\end{bmatrix} .
\end{aligned}
\end{equation}
We only take the real part of $H$ and $\bfJ$ and the imaginary part of $\bfE$ and $q$ into consideration. In order to have this analytical solution \eqref{eq:anaSol}, one needs to set the length of the square $L=\frac{\sqrt{2}\pi}{k}$ and $\beta^2=\frac{\omega^2}{k^2}$. The convergence history of the HDG method with interpolation order $\mathbb{P}_p \ (p=1,2,3)$ is given in Table \ref{tbl:convRect} and Figure \ref{fig:convRect}. Mesh size $h$ is the edge length of elements associated to the boundary $\partial\Omega_{\Box}$. The convergence orders are calculated by
$$
\frac{\log(\|\bfE^a-\bfE^h\|_{\Omega_{\Box}}^{h_2}/\|\bfE^a-\bfE^h\|_{\Omega_{\Box}}^{h_1})}{\log(h_2/h_1)},
$$
where $h_1$ and $h_2$ denote a coarse and a refined mesh size, respectively. From Table~\ref{tbl:convRect} and Figure~\ref{fig:convRect}, we observe that the proposed HDG method has an optimal convergence order which is $p+1$ for HDG-$\mathbb{P}_p$.

\begin{table}
\centering \caption{Convergence results for the cavity problem.}
\begin{tabular}{ccccccc}\label{tbl:convRect}\\
\hline
&  \multicolumn{2}{c}{HDG-$\mathbb{P}_1$} & \multicolumn{2}{c}{HDG-$\mathbb{P}_2$} & \multicolumn{2}{c}{HDG-$\mathbb{P}_3$} \\
\hline
$h$ & error & order & error & order & error & order\\
\hline
$5\time10^{-8}$ & $1.67\times10^{-9}$ & - & $4.52\times10^{-10}$ & - & $2.04\times10^{-11}$ & - \\
$2.5\time10^{-8}$ & $4.10\times10^{-10}$ & 2.0 & $5.61\times10^{-11}$ & 3.0 & $1.28\times10^{-12}$ & 4.0 \\
$1.25\time10^{-8}$ & $9.98\times10^{-11}$ & 2.0 & $7.52\times10^{-12}$ & 3.0 & $7.78\times10^{-14}$ & 4.0 \\
$6.25\time10^{-9}$ & $2.40\times10^{-12}$ & 2.1 & $9.11\times10^{-13}$ & 3.0 & $5.03\times10^{-15}$ & 4.0 \\
\hline
\end{tabular}
\end{table}

\subsection{Benchmark problem: a cylindrical plasmonic nanowire}
As our benchmark problem we consider the plasmonic behavior of a cylindrical nanowire. This has been used as a convenient benchmark problem for other numerical methods before~\cite{ToscanoEtal2012,HiremathEtal2012} because  analytical solutions exist both for the local and for the NHD models, see the derivation in Ref.~\cite{Ruppin2001}. 
We make use of the fact that the analytical Mie solution of Ref.~\cite{Ruppin2001} allows making the nonlocal parameter $\beta$ complex-valued. This enables us to benchmark our HDG simulations against exact analytical results for the GNOR model as well. (For comparison, optical properties of a sphere in the GNOR model, based on exact Mie results, are discussed in Ref.~\cite{RazaEtal2015JPCM}.)

For the NHD model, the configuration of the nanowire is taken to be the same as that in the first test in~\cite{HiremathEtal2012}: the radius of the cylinder is 2\,nm, no interband transitions are considered, the plasma frequency $\omega_p=8.65\times10^{15}$, the damping constant $\gamma=0.01\omega_p$, the Fermi velocity $v_F=1.07\times 10^{6}$, and $\beta^2=\frac{3}{5}v_F^2$. For the GNOR model, we use the same parameters and furthermore we take $D=2.04\times10^{-4}$ ~\cite{MortensenGNOR}. An artificial absorbing boundary is set to be a concentric circle with a radius of $100$\,nm. 

As our benchmark observable we will calculate the Extinction Cross Section (ECS, $\sigma_{\text{ext}}$), which is given by the sum of the scattering cross section $\sigma_\text{sca}$ and the absorption cross section $\sigma_\text{abs}$~\cite{BergEtal2009},
$$
\sigma_\text{ext} = \sigma_\text{sca} + \sigma_\text{abs}.
$$
More precisely, for the cylindrical nanowire we consider the  extinction cross section per wire length, which actually has the units of a length. We scale this quantity by the diameter $2 r$ of the nanowire to obtain a dimensionless normalized extinction cross section that we denote by $\sigma_\text{ext}$. It can be expressed as the sum of scaled scattering and absorption cross sections,
$$
\sigma_\text{sca} = \frac{1}{2r}Re \oint_S(\bfE^\text{sca}\times \overline{H}^\text{sca})\cdot\bfn \, dS, \ \text{ and } \ \sigma_\text{abs} = -\frac{1}{2r}Re \oint_S(\bfE^\text{tot}\times \overline{H}^\text{tot})\cdot\bfn \, dS. 
$$
Here the integrations are performed along a closed path around the nanowire, and $Re$ denotes the real part.

\begin{figure}
\centering
\subfigure[]{\includegraphics[width=0.45\textwidth,height=0.25\textheight]{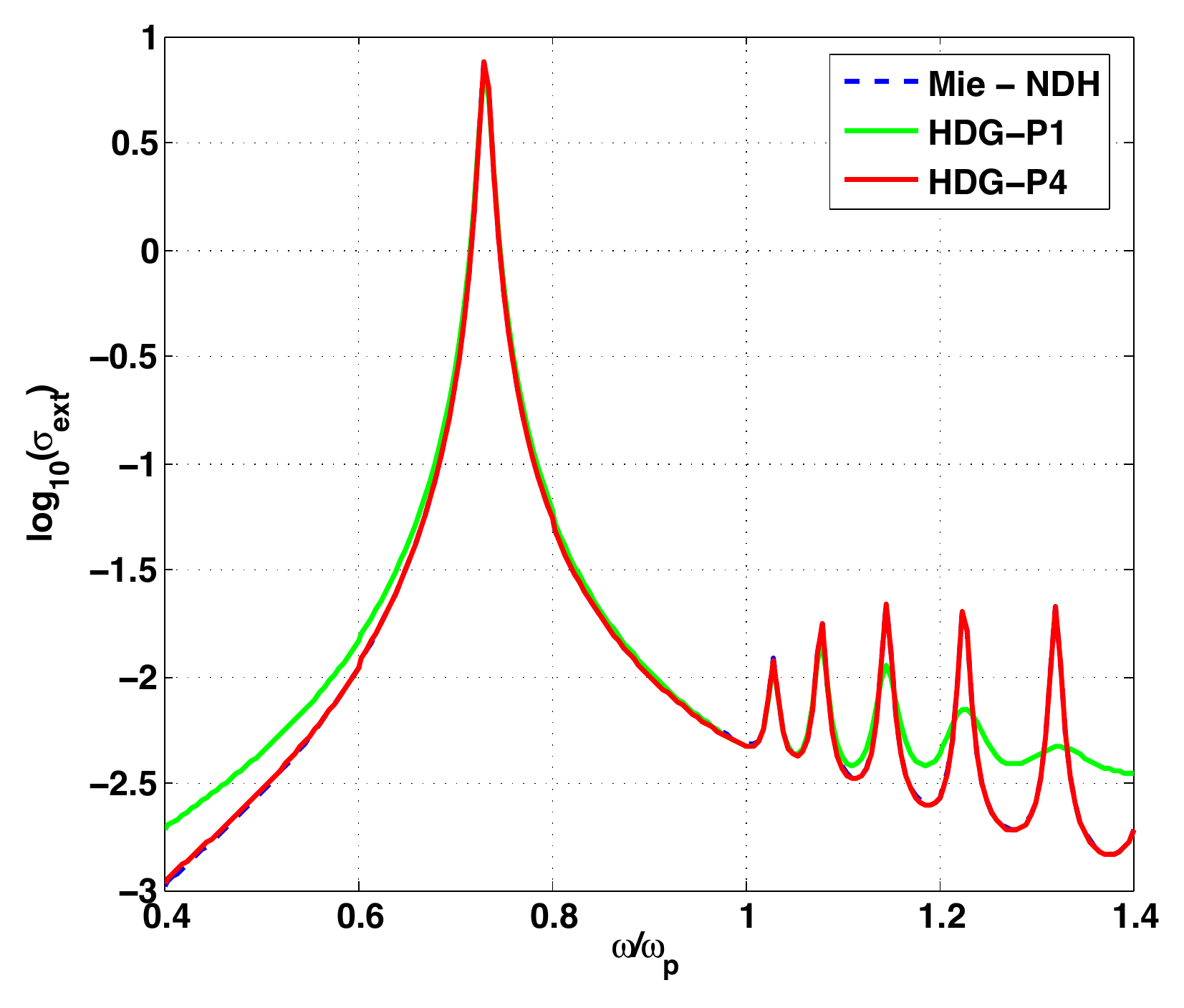}}
\subfigure[]{\includegraphics[width=0.45\textwidth,height=0.25\textheight]{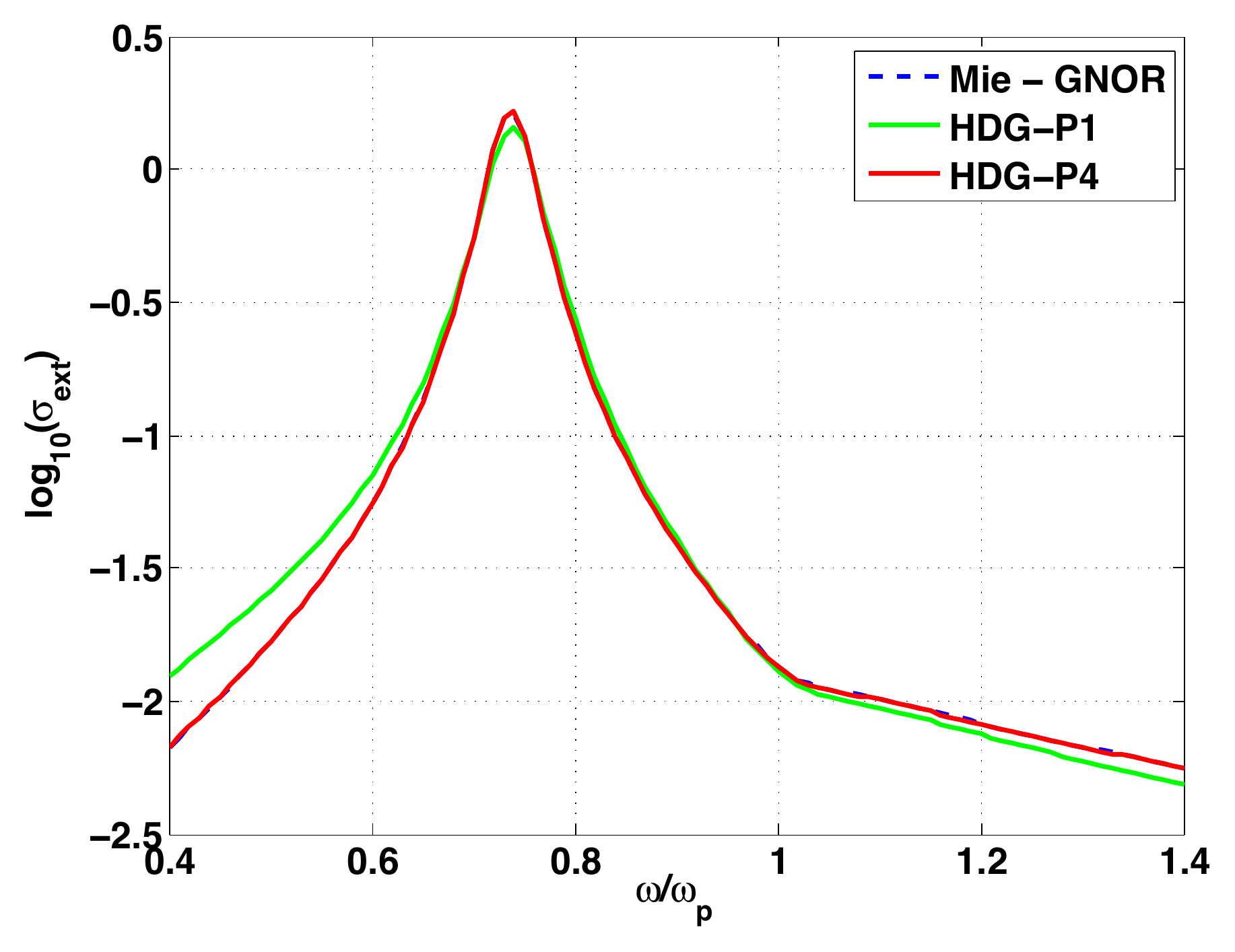}}
\caption{Extinction cross section of a Na cylinder with 2 nanometer radius in a free-space background, for a TM-polarized normally incident plane wave. The cylinder is described both by (a) the NHD model and by (b) the GNOR model. The simple wire geometry serves as an excellent benchmark problem: Analytically exact calculations (nonlocal Mie theory) are compared with  HDG methods of different interpolation order.}\label{fig:Cext_nanowire}
\end{figure}

\begin{figure}
\centering
\includegraphics[width=0.95\textwidth]{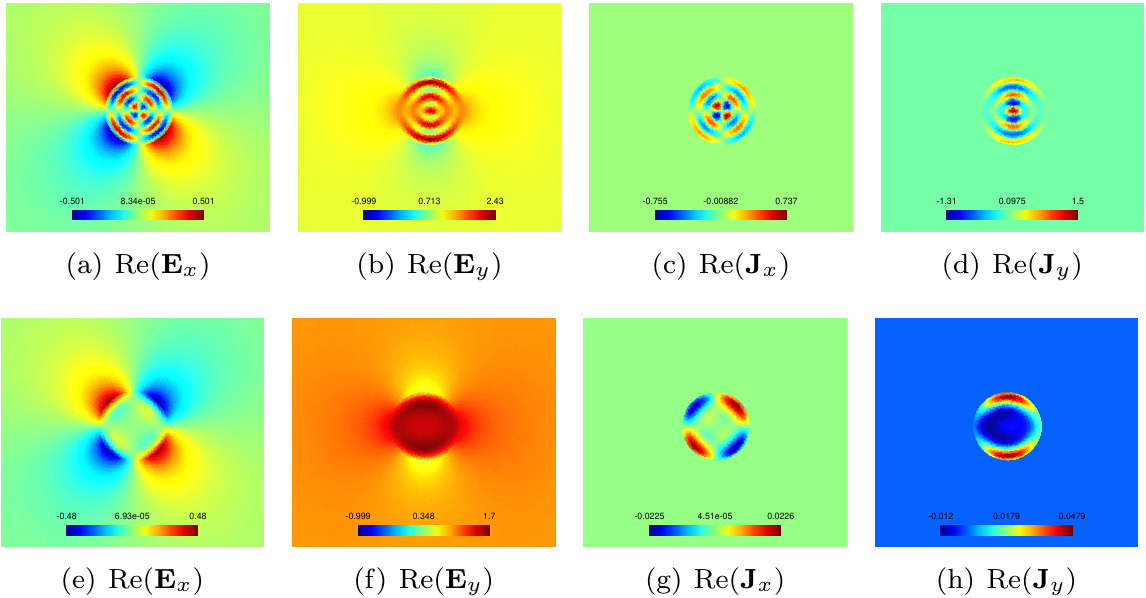} 
\caption{The electric-field and current-density distributions of the light-matter interaction of a Na nanowire. In the upper row, we show the distributions on the fourth-order nonlocal resonance at $\omega/\omega_p=1.227$ for the NHD model. For comparison, in the bottom row we show the corresponding distributions for the GNOR model.}\label{fig:Fields_nanowire}
\end{figure}

The simulations are performed on a mesh with 4,513 nodes, 8,896 elements and 13,280 edges of which 722 edges are located inside the nanostructure. The ECS is presented in Figure~\ref{fig:Cext_nanowire}. Curvilinear treatment is employed for high-order accuracy, where the curved edges are geometrically approximated by second-order curves instead of straight lines~\cite{LLP2013}. From Figure~\ref{fig:Cext_nanowire} we can observe that the fourth-order HDG method produces an ECS curve that matches the analytical solution very well. By contrast, the first-order method is not accurate enough on this mesh. Contour plots of the electric field and the current density are presented in Figure~\ref{fig:Fields_nanowire}. These results match well with corresponding results in Ref.~\cite{HiremathEtal2012} despite the lower resolution, probably because our simulation is performed on a coarser mesh. Comparing the two subfigures in Figure~\ref{fig:Cext_nanowire}, we also find that the ECS curve for the GNOR model is smoother than for the NHD model. But this has a physical rather than a numerical origin. In particular  the standing bulk plasmon resonances  above the plasma frequency in the NHD model are essentially washed out by the introduced diffusion in the GNOR model. The ECS curves of HDG-$\mathbb{P}_2$ and HDG-$\mathbb{P}_3$ are are not presented in Figure~\ref{fig:Cext_nanowire}, but we found that they they lie in between the displayed curves of  HDG-$\mathbb{P}_1$ and  HDG-$\mathbb{P}_4$. 

In our 2D simulations, we use a sparse direct solver MUMPS~\cite{AmestoyEtal2000} to solve the resulting systems of linear equations. We need to solve a linear system at each frequency. The computational performance mainly relies on the size of the coefficient matrices, \emph{i.e.} the number of degrees of freedom (\#DOF). The computational performance for one frequency is given in Table~\ref{tbl:perf}, where $t_\text{construction}$ denotes the CPU time for construction the matrices, $t_\text{factorization}$ denotes the CPU time used by MUMPS for the factorization of the coefficient matrix $A$ \eqref{eq:linSys}, and memory denotes the memory consumed by MUMPS. From Table~\ref{tbl:perf} we can see that the HDG-$\mathbb{P}_4$ is more expensive than HDG-$\mathbb{P}_1$ in CPU time for both construction and factorization. However, high-order methods are preferable because they costs less for the same accuracy~\cite{LLP2013}.

\begin{table}
\centering \caption{Computational performance of the nanowire problem.}
\begin{tabular}{ccccc}\label{tbl:perf}\\
\hline
  & \#DOF & $t_\text{construction}$ (second) & $t_\text{factorization}$  (second) & memory (MB) \\
HDG-$\mathbb{P}_1$ & 28,260 & 0.067 & 0.36 & 74 \\
HDG-$\mathbb{P}_4$ & 70,650 & 2.4 & 3.3  & 418\\
\hline
\end{tabular}
\end{table}

\subsection{Dimer of cylindrical nanowires}
\label{sec:dimer}
\begin{figure}[htbp]
  \centering
  \subfigure[Configuration]{\includegraphics[width=0.45\textwidth]{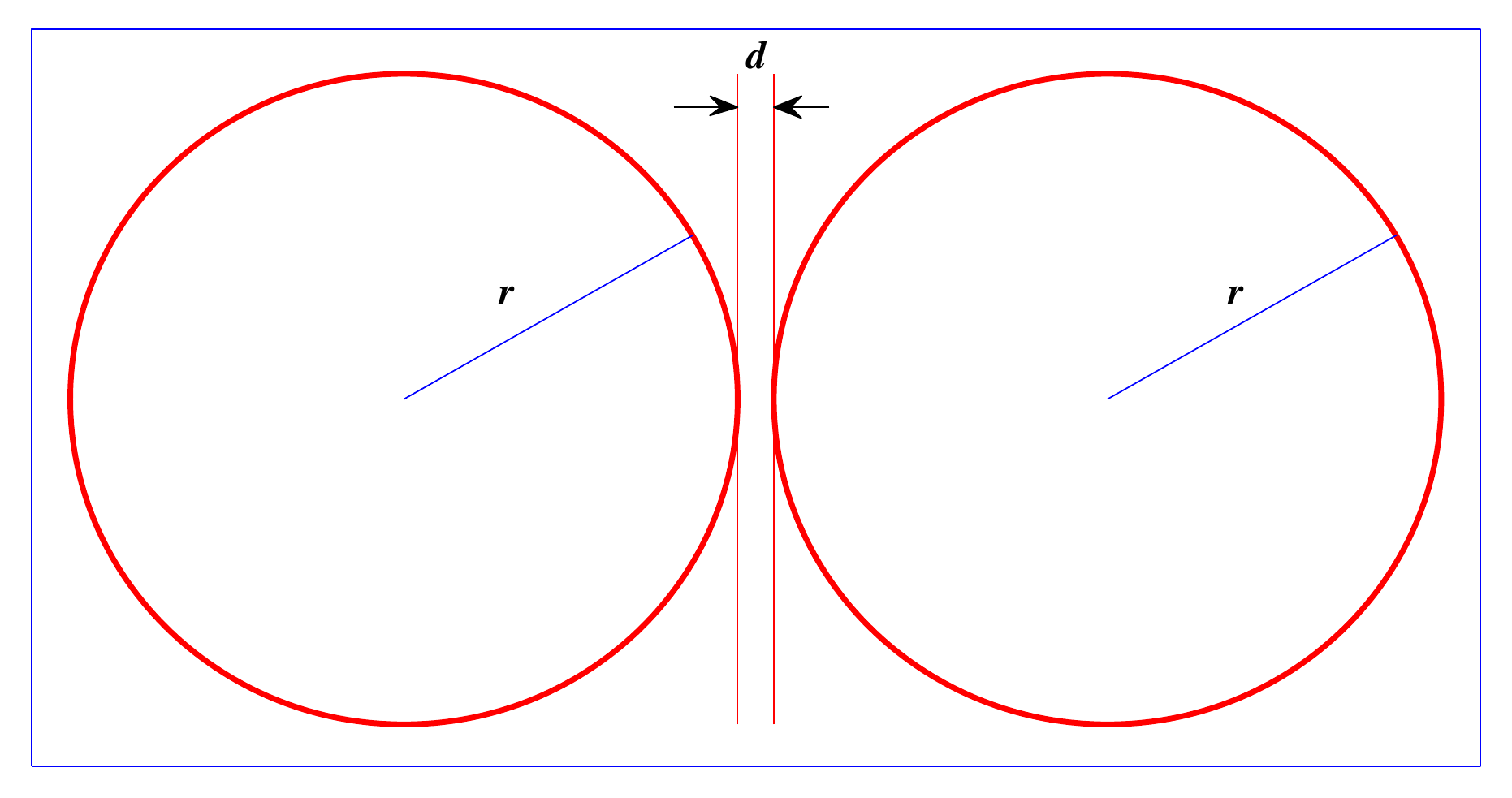}}
  \subfigure[A mesh]{\includegraphics[width=0.35\textwidth]{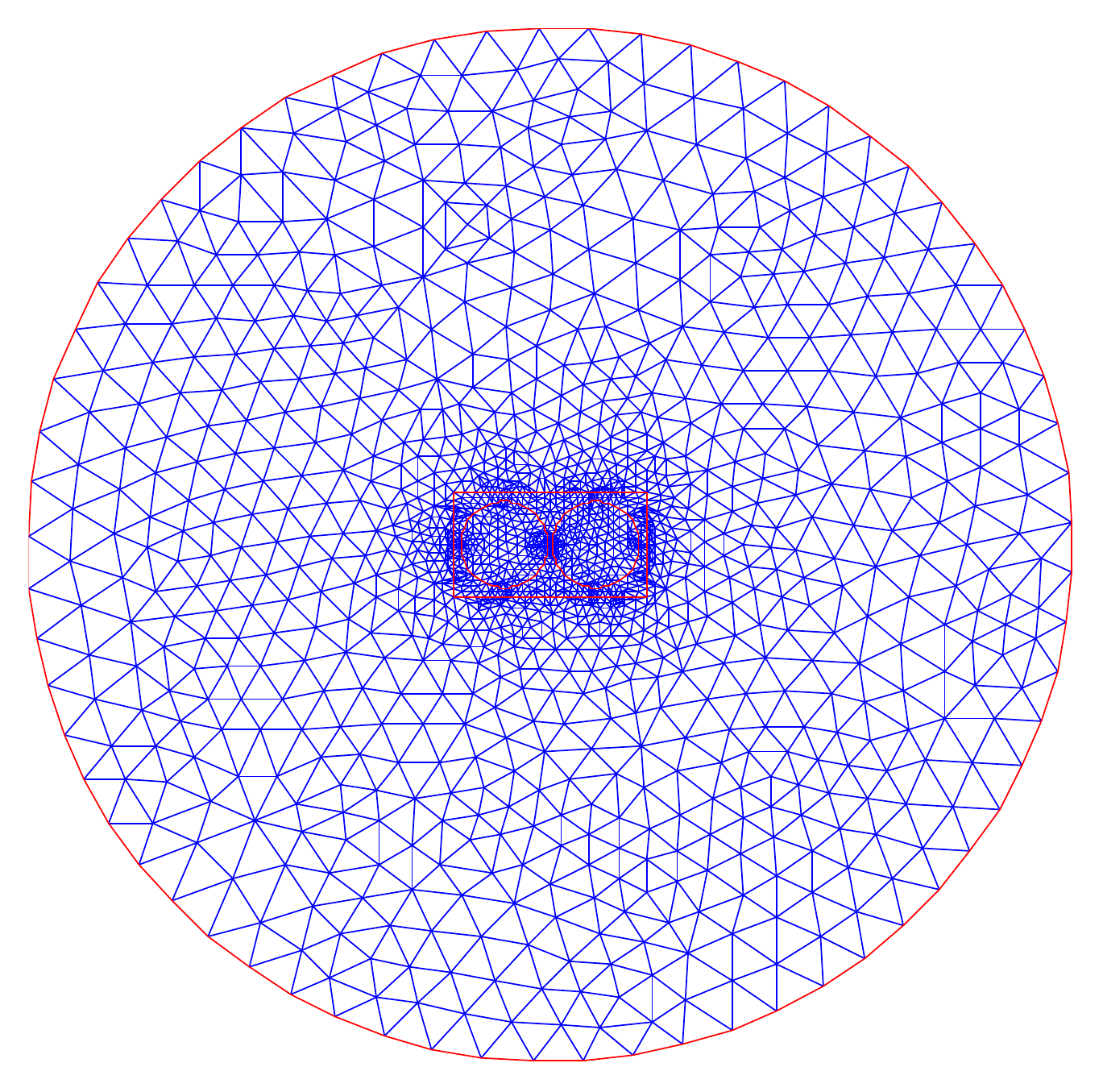}}
  \caption{A cylindrical gold dimer nanowire. Panel (a): The geometry of the  nanowire dimer.  An imaginary rectangle  is introduced around the dimer for the  calculation of the cross section.   Both nanowires have a radius $r=30$\,nm and their gap distance is $d=3$\,nm. Panel (b): A typical mesh. The  large circle is the artificial absorbing boundary.}
  \label{fig:dimerConfig}
\end{figure}

Plasmonic dimer structures with small gaps are both experimentally interesting and computationally challenging because of high field enhancements in the gap region~\cite{ToscanoEtal2012,Gallinet2015,RazaEtal2015JPCM}. Here we present our HDG simulations of a cylindrical gold dimer geometry as shown in Figure \ref{fig:dimerConfig}(a), and this particular configuration is from Ref.~\cite{RazaEtal2015}. A typical mesh is shown in Figure~\ref{fig:dimerConfig}(b).
On a mesh with 5,829 nodes, 11,520 triangles and 17,348 edges with 3,712 edges inside the nanostructure, we calculate the ECS curve by  HDG-$\mathbb{P}_4$. The size of matrix for HDG-$P1$ is $105,300\times105,300$, the matrix construction CPU time is 5.2 seconds, the factorization CPU time is 6.9 seconds for one frequency, and the memory cost is 717 MB.

\begin{figure}
\centering
\includegraphics[width=0.45\textwidth]{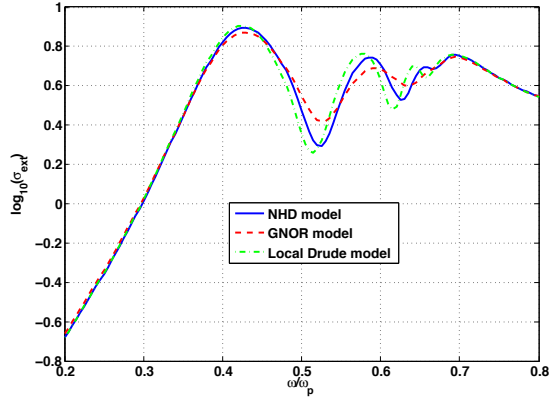}
\caption{Comparison of extinction cross sections of a gold dimer as calculated with the local Drude model, the NHD model and the GNOR model. The configuration is shown in Figure~\ref{fig:dimerConfig} and the material parameters are given in subsection~\ref{sec:dimer}. }\label{fig:dimECS}
\end{figure}

\begin{figure}
\centering
\includegraphics[width=0.95\textwidth]{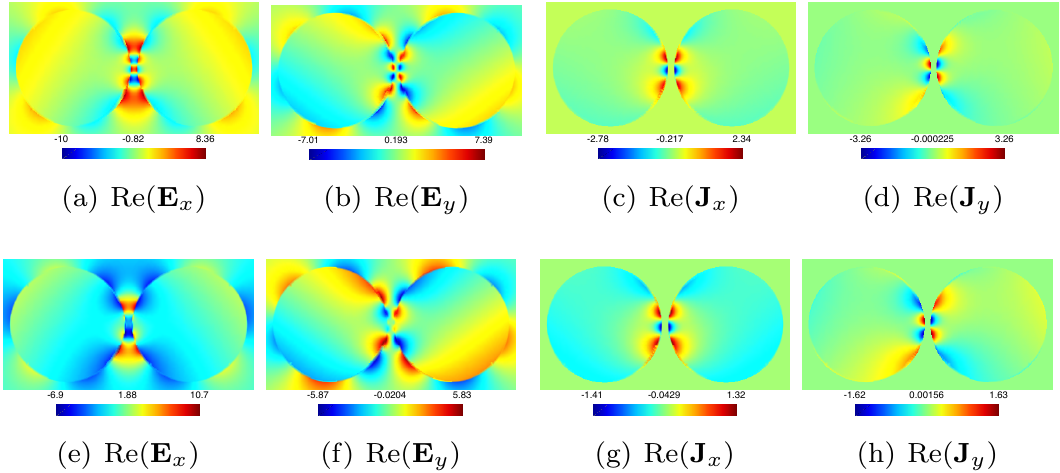}
\caption{Various field distributions in the gold dimer when illuminated by a plane wave of light. On the top line we show the distributions at the third SPR  of the NHD model, at $\omega/\omega_p=0.66$. On the second line we show the corresponding distributions in the GNOR model at the same frequency.}\label{fig:dimerFields}
\end{figure}

For the material properties gold we use the same values as in Ref.~\cite{ToscanoEtal2012}: the plasma frequency $\omega_p=1.34\times10^{16}$, damping constant $\gamma=1.14\times 10^{14}$, the Fermi velocity $v_F=1.39\times 10^{6}$, and the nonlocal parameter $\beta$ is determined by $\beta^2=\frac{3}{5}v_F^2$. The incoming plane wave  of light is incident perpendicular to the line connecting the centers of the two circles, with a linear polarization parallel to this line (TM- or $p$-polarization). A comparison of the ECS curves is presented in Figure~\ref{fig:dimECS}. Overall, there are small but clear differences, illustrating that nonlocal response effects occur even for dimer structures for which the corresponding monomers ($r= 30$\,nm nanowires) would show essentially no nonlocal effects~\cite{ToscanoEtal2012}. Both nonlocal models have blueshifted resonances as compared to the local model, and resonances in the GNOR model are less pronounced than in the local and NHD models. For smaller gap sizes, nonlocal blueshifts are larger and resonances are broadened more (the latter only in the GNOR model).
Field distributions at the same particular frequency for the NHD  and the GNOR models are compared in Figure~\ref{fig:dimerFields}. The figure illustrates the generic features that the GNOR model washes out some finer details of the field distributions, and also that minimal and maximal field values lie closer together in the GNOR model.

\section{Conclusions}
\label{sec:con}
This paper introduces a HDG method to solve the nonlocal hydrodynamic Drude model and the GNOR model, both of which are  often employed to describe light-matter interactions of nanostructures. The numerical fluxes are expressed in terms of two newly introduced hybrid terms. Only the hybrid unknowns are involved in the global problem. The local problems are solved element-by-element once the hybrid terms are obtained. The proposed HDG formulations naturally couple the hard-wall boundary condition.  Numerical results indicate that the HDG method converges at the optimal rate. Our benchmark simulations for a cylindrical nanowire and our calculations for a dimer structure show that the HDG method is a promising method in nanophotonics. Building on these results, in the near future we plan to generalize our computations to 3D structures, and to introduce domain decomposition and model order reduction into nanophotonic computations. 

\section*{Acknowledgments}
The first author was supported by the NSFC (11301057) and the Fundamental Research Funds for the Central Universities (ZYGX2014J082). 
N. A. M. and M. W. acknowledge support from the Danish Council for Independent Research (FNU 1323-00087). 
M. W. acknowledges support from the Villum Foundation via the VKR Centre of
Excellence NATEC-II.
The Center for Nanostructured Graphene is sponsored by
the Danish National Research Foundation, Project DNRF103.

%---------------------------------------------------------------------

\end{document}